\documentclass[twoside,leqno,10pt, A4]{amsart}
\usepackage{amsfonts}
\usepackage{amsmath}
\usepackage{amscd}
\usepackage{amssymb}
\usepackage{amsthm}
\usepackage{amsrefs}
\usepackage{latexsym}
\usepackage{mathrsfs}
\usepackage{bbm}
\usepackage{enumerate}
\usepackage{graphicx}

\usepackage{amsfonts}
\usepackage{amsmath}
\usepackage{amscd}
\usepackage{amssymb}
\usepackage{amsthm}
\usepackage{amsrefs}
\usepackage{latexsym}
\usepackage{mathrsfs}
\usepackage{bbm}
\usepackage{amscd}
\usepackage{amssymb}
\usepackage{amsthm}
\usepackage{amsrefs}
\usepackage{latexsym}
\usepackage{mathrsfs}
\usepackage{bbm}
\usepackage{enumerate}
\usepackage{graphicx}
\usepackage{color}
\begin{document}

\newtheorem{theorem}[subsection]{Theorem}
\newtheorem{proposition}[subsection]{Proposition}
\newtheorem{lemma}[subsection]{Lemma}
\newtheorem{corollary}[subsection]{Corollary}
\newtheorem{conjecture}[subsection]{Conjecture}
\newtheorem{prop}[subsection]{Proposition}
\numberwithin{equation}{section}
\newcommand{\mr}{\ensuremath{\mathbb R}}
\newcommand{\mc}{\ensuremath{\mathbb C}}
\newcommand{\dif}{\mathrm{d}}
\newcommand{\intz}{\mathbb{Z}}
\newcommand{\ratq}{\mathbb{Q}}
\newcommand{\natn}{\mathbb{N}}
\newcommand{\comc}{\mathbb{C}}
\newcommand{\rear}{\mathbb{R}}
\newcommand{\prip}{\mathbb{P}}
\newcommand{\uph}{\mathbb{H}}
\newcommand{\fief}{\mathbb{F}}
\newcommand{\majorarc}{\mathfrak{M}}
\newcommand{\minorarc}{\mathfrak{m}}
\newcommand{\sings}{\mathfrak{S}}
\newcommand{\fA}{\ensuremath{\mathfrak A}}
\newcommand{\mn}{\ensuremath{\mathbb N}}
\newcommand{\mq}{\ensuremath{\mathbb Q}}
\newcommand{\half}{\tfrac{1}{2}}
\newcommand{\f}{f\times \chi}
\newcommand{\summ}{\mathop{{\sum}^{\star}}}
\newcommand{\chiq}{\chi \bmod q}
\newcommand{\chidb}{\chi \bmod db}
\newcommand{\chid}{\chi \bmod d}
\newcommand{\sym}{\text{sym}^2}
\newcommand{\hhalf}{\tfrac{1}{2}}
\newcommand{\sumstar}{\sideset{}{^*}\sum}
\newcommand{\sumprime}{\sideset{}{'}\sum}
\newcommand{\sumprimeprime}{\sideset{}{''}\sum}
\newcommand{\sumflat}{\sideset{}{^\flat}\sum}
\newcommand{\shortmod}{\ensuremath{\negthickspace \negthickspace \negthickspace \pmod}}
\newcommand{\V}{V\left(\frac{nm}{q^2}\right)}
\newcommand{\sumi}{\mathop{{\sum}^{\dagger}}}
\newcommand{\mz}{\ensuremath{\mathbb Z}}
\newcommand{\leg}[2]{\left(\frac{#1}{#2}\right)}
\newcommand{\muK}{\mu_{\omega}}
\newcommand{\thalf}{\tfrac12}
\newcommand{\lp}{\left(}
\newcommand{\rp}{\right)}
\newcommand{\Lam}{\Lambda_{[i]}}
\newcommand{\lam}{\lambda}
\def\L{\fracwithdelims}
\def\om{\omega}
\def\pbar{\overline{\psi}}
\def\phi{\varphi}
\def\lam{\lambda}
\def\lbar{\overline{\lambda}}
\newcommand\Sum{\Cal S}
\def\Lam{\Lambda}
\newcommand{\sumtt}{\underset{(d,2)=1}{{\sum}^*}}
\newcommand{\sumt}{\underset{(d,2)=1}{\sum \nolimits^{*}} \widetilde w\left( \frac dX \right) }

\theoremstyle{plain}
\newtheorem{conj}{Conjecture}
\newtheorem{remark}[subsection]{Remark}

\providecommand{\re}{\mathop{\rm Re}}
\providecommand{\im}{\mathop{\rm Im}}
\def\cI{\mathcal{I}}
\def\cL{\mathcal{L}}
\def\E{\mathbb{E}}

\makeatletter
\def\widebreve{\mathpalette\wide@breve}
\def\wide@breve#1#2{\sbox\z@{$#1#2$}%
     \mathop{\vbox{\m@th\ialign{##\crcr
\kern0.08em\brevefill#1{0.8\wd\z@}\crcr\noalign{\nointerlineskip}%
                    $\hss#1#2\hss$\crcr}}}\limits}
\def\brevefill#1#2{$\m@th\sbox\tw@{$#1($}%
  \hss\resizebox{#2}{\wd\tw@}{\rotatebox[origin=c]{90}{\upshape(}}\hss$}
\makeatletter

\title[Lower bounds for negative moments of quadratic Dirichlet $L$-functions]{Lower bounds for negative moments of quadratic Dirichlet $L$-functions}

%%\date{\today}
\author{Peng Gao}
\address{School of Mathematical Sciences, Beihang University, Beijing 100191, P. R. China}
\email{penggao@buaa.edu.cn}

\begin{abstract}
 We establish lower bounds for the $2k$-th moment of families of quadratic Dirichlet $L$-functions at the central point for all real $k<0$, assuming a conjecture of S. Chowla on the non-vanishing of these $L$-values.
\end{abstract}

\maketitle



\section{Introduction}
\label{sec 1}

   Central values of $L$-functions have important arithmetic implications and their non-vanishing issues have received much attention. For quadratic Dirichlet $L$-functions, it is conjectured by S. Chowla \cite{chow} that $L(1/2, \chi) \neq 0$ for every primitive real Dirichlet character $\chi$. Throughout paper, we shall refer to this conjecture as Chowla's Conjecture.

  There has been considerable progress made towards Chowla's Conjecture. M. Jutila \cite{Jutila} initiated the study on moments of families of quadratic Dirichlet $L$-functions to show that Chowla's Conjecture is true for infinitely many such $L$-functions.
A well-known result of K. Soundararajan \cite{sound1} using mollified moments shows that $L(1/2,\chi_{8d})\neq 0$ for at least $87.5\%$ of the real characters $\chi_{8d}$ for $d$ odd and square-free. Here $\chi_{8d} = \left(\frac{8d}{\cdot} \right)$ is the Kronecker symbol. In \cite{B&P}, S. Baluyot and K. Pratt proved that more than $9 \%$ of the members of the quadratic family of Dirichlet $L$-functions with prime moduli do not vanish at the central point.

 On the other hand, the density conjecture of N. Katz and P. Sarnak \cites{KS1, K&S} would imply that $L(1/2, \chi) \neq 0$ for almost all quadratic Dirichlet character $\chi$. In this direction, A. E. \"{O}zluk and C. Snyder  \cite{O&S} computed the one level density for a corresponding family to show that $L(1/2, \chi_d) \neq 0$ for at least $15/16$ of the fundamental discriminants $|d| \leq X$. Their result assumes the generalized Riemann hypothesis (GRH). Using an optimal test function as in \cites{B&F, ILS}, this percentage can be improved to be  $(19-\cot \frac{1}{4})/16$. Also via computing the one level density, it is shown by J. C. Andrade and S. Baluyot \cite[Theorem 3]{A&B} that at least $75\%$ of the family of quadratic Dirichlet $L$-functions with prime moduli does not vanish at the central point under GRH.

  As moments of families of $L$-functions can be applied effectively to address the non-vanishing issues of $L$-functions, the study on moments of  $L$-functions becomes an important topic in analytical number theory.  There are now precise conjectures concerning the asymptotic behaviors of non-negative moments, owing to the work of J. P. Keating and N. C. Snaith \cite{Keating-Snaith02} using random matrix theory, and also to the work of J. B. Conrey, D. W. Farmer, J. P. Keating, M. O. Rubinstein and N. C. Snaith \cite{CFKRS}. Moreover, systematic ways towards establishing sharp lower bounds concerning these moments are developed by Z. Rudnick and K. Soundararajan \cites{R&Sound, R&Sound1}, M. Radziwi{\l\l} and K. Soundararajan \cite{Radziwill&Sound1}, W. Heap and K. Soundararajan \cite{H&Sound}. Also, corresponding methods for establishing sharp upper bounds concerning these moments are developed by K. Soundararajan \cite{Sound2009} with a refinement by A. J. Harper \cite{Harper} (under GRH), M. Radziwi{\l\l} and K. Soundararajan \cite{Radziwill&Sound}.

  For the family of quadratic Dirichlet $L$-functions, M. Jutila \cite{Jutila} obtained asymptotic formulas for the first and the second moments for the family with general moduli and for the first moment for the family with prime moduli. Third and fourth moment for the same family with general moduli were evaluated asymptotically by K. Soundararajan \cite{sound1}  unconditionally and under GRH by Q. Shen \cite{Shen}, respectively.  Subsequent improvements on the error terms concerning these moments can be found in \cites{ViTa, DoHo, Young1, Young2, sound1, Sono}.  An asymptotic formula for the second moment with prime moduli was obtained by S. Baluyot and K. Pratt \cite{B&P} under GRH.

Moreover, sharp bounds for the $k$-th moment of the above families have
been extensively studied in  \cites{Harper, Sound2009,
R&Sound, R&Sound1, Radziwill&Sound1, Gao2021-2, Gao2021-3}. For the family with general moduli, sharp lower bounds (resp. upper bounds) are established for all real $k \geq 0$ (resp. $0 \leq k \leq 2$) unconditionally, while sharp upper bounds are obtained for all
real $k > 2$ under GRH. For the family with prime moduli or prime-related moduli,  sharp upper and lower bounds were obtained by S. Baluyot and K. Pratt \cite{B&P} for the second unconditionally and third moment under certain assumptions, by P. Gao and L. Zhao \cite{G&Zhao11} for all $k$-th moment with $k \geq 0$  under GRH.

   In this paper, we are interested in the negative moments of the above families of quadratic Dirichlet $L$-functions. There are relatively fewer
results concerning these negative moments compared to those on non-negative moments, partially due to the reason that the $L$-functions may vanish at the central point. Also, the behaviour of the negative moments may be more difficult to predict compared to that of the positive ones. Take for example the $2k$-th discrete moment of the derivative of the Riemann zeta function $\zeta(s)$ at nontrivial zeros given by
\begin{align*}
J_k(T) =\frac{1}{N(T)}\sum_{0<\Im(\rho)\le T}|\zeta'(\rho)|^{2k},
\end{align*}
  where we denote $\rho$ for the nontrivial zeros of $\zeta(s)$ and $N(T) =\sum_{0<\Im(\rho)\le T}1$.

 A conjecture due to S. M. Gonek \cite{Gonek1} and D. Hejhal \cite{Hejhal} independently asserts that for any real $k$,
\begin{align*}
%%\label{Jksim}
J_k(T) \asymp (\log T)^{k(k+2)}.
\end{align*}

  In \cite{HKO}, C. P. Hughes, J. P. Keating and N. O'Connell gave a more precise prediction concerning $J_k(T)$ using the random matrix theory to show that the order of magnitude undergoes a phase change at $k=-3/2$. This prediction was further confirmed by H. M. Bui, S. M. Gonek and M. B. Milinovich \cite{BGM} using the hybrid Euler-Hadamard product. In \cite{HLZ}, W. Heap, J. Li and J. Zhao proved that $J_k(T) \gg (\log T)^{k(k+2)}$ for all rational $k<0$ assuming the Riemann hypothesis (RH) and all zeros of $\zeta(s)$ are simple.

  For the case of quadratic Dirichlet $L$–functions, random matrix theory computations due to P. J. Forrester and J. P. Keating \cite{FK} also suggest certain phase changes in the asymptotic formulas for the $2k$-th moment when $2k = -(2j + 1/2)$ for any positive integer $j$. In \cite{BFK21, Florea21}, upper bounds for negative moments of quadratic Dirichlet $L$-functions over function fields were obtained at points that are slightly shifted away from the central point.

  It is the aim in our paper to establish lower bounds for negative moments of families of quadratic Dirichlet $L$–functions. For this, it is natural to assume Chowla's Conjecture, so that negative powers of central values of these $L$-functions become meaningful. We shall consider families with general moduli as well as with prime-related moduli. More precisely, we are interested in the following families of $L$-functions at the central point:
\begin{align*}
  \big\{ L(\half,\chi_{8d}) : d \ \text{odd and square-free} \big\} \ \ \text{and} \ \ \big\{ L(\half,\chi_{8p}) : p \ \text{odd prime} \big\}.
\end{align*}
  Here we note that $\chi_{8d}$ is a primitive Dirichlet character (see \cite[\S 5]{Da}) for any odd, square-free $d$. We shall further reserve the letter $p$ for a prime throughout the paper.

 Our result is as follows.
\begin{theorem}
\label{thmlowerbound}
   Assume the truth of Chowla's Conjecture. Let $k$ be any negative real number and $X$ be a large real number. We have
\begin{align*}
   \sumstar_{\substack{ 0<d \leq X \\ (d,2)=1}}|L(\tfrac{1}{2},\chi_{8d})|^{2k} \gg_k X(\log X)^{\frac{2k(2k+1)}{2}},
\end{align*}
  where we denote $\sumstar$ for the sum over square-free integers.

  Assume further the truth of RH for $\zeta(s)$ and GRH for $L(s, \chi_{8p})$ for all odd primes $p$, then we have
\begin{align*}
%%\label{lowerbound}
\begin{split}
    \sum_{\substack{2< p \leq X}} ( \log p) |L(\tfrac{1}{2},\chi_{8p})|^{2k}  \gg_k & X(\log X)^{\frac{2k(2k+1)}{2}}.
\end{split}
\end{align*}
\end{theorem}

   Our proof of Theorem \ref{thmlowerbound} is based on a variant of the lower bounds principle of W. Heap and K. Soundararajan \cite{H&Sound} and is inspired by the work of W. Heap, J. Li and J. Zhao \cite{HLZ} on lower bounds of $J_k(T)$ for negative $k$. The proof also utilizes much of the work done in \cite{Gao2021-3}. In view of the prediction from random matrix theory given in \cite{FK}, it is plausible that the bounds given in Theorem \ref{thmlowerbound} are sharp for $-5/2 \leq 2k <0$, provided that Chowla's Conjecture is true.

\section{Proof of Theorem \ref{thmlowerbound}}
\label{sec 2'}

    We follow the treatment in \cite{Gao2021-3} to denote $\Phi$ for a smooth function, compactly supported on $[1/8, 7/8]$ such that $0 \leq \Phi(x) \leq 1$ for all $x$ and that  $\Phi(x) = 1$ for $x \in [1/4, 3/4]$. Upon dividing $d$ or $p$ into dyadic blocks, we see that in order to prove Theorem \ref{thmlowerbound}, it suffices to show that
\begin{align}
\label{orderofmagdag}
   \sumstar_{\substack{(d,2)=1}}|L(\tfrac{1}{2},\chi_{8d})|^{2k}\Phi(\frac dX) \gg_k &  X(\log X)^{\frac{2k(2k+1)}{2}}, \\
\label{orderofmagp}
   \sum_{\substack{(p,2)=1}}( \log p)|L(\tfrac{1}{2},\chi_{8p})|^{2k}\Phi(\frac pX) \gg_k & X(\log X)^{\frac{2k(2k+1)}{2}}.
\end{align}

    Denote $N, M$ for two large natural numbers depending on $k$ only and let $\{ \ell_j \}_{1 \leq j \leq R}$ be a sequence of even natural numbers such that $\ell_1= 2\lceil N \log \log X\rceil$ and $\ell_{j+1} = 2 \lceil N \log \ell_j \rceil$ for $j \geq 1$, where we set $R$ to be the largest natural number satisfying $\ell_R >10^M$.  We may also choose $M$ large enough to ensure $\ell_{j} > \ell_{j+1}^2$ for all $1 \leq j \leq R-1$. This also implies that
\begin{align*}
%%\label{sumoverell}
  R \ll \log \log \ell_1, \quad  \sum^R_{j=1}\frac 1{\ell_j} \leq \frac 2{\ell_R}.
\end{align*}

   Let ${ P}_1$ be the set of odd primes not exceeding $X^{1/\ell_1^2}$,
${ P_j}, 2\le j\le R$ be the set of primes lying in the interval $(X^{1/\ell_{j-1}^2}, X^{1/\ell_j^2}]$. For any odd and square-free $d$, we denote
\begin{equation*}
%% \label{defP}
{\mathcal P}_j(d) = \sum_{\substack{q\in P_j}} \frac{1}{\sqrt{q}} \chi_{8d}(q).
\end{equation*}

  We set for any integer $\ell \geq 0$ and any $x \in \mr$,
\begin{equation*}
%%\label{E_ell}
E_{\ell}(x) = \sum_{j=0}^{\ell} \frac{x^{j}}{j!}.
\end{equation*}

  Then we define for any $\alpha \in \mr$,
\begin{align}
\label{defN}
{\mathcal N}_j(d, \alpha) = E_{\ell_j} (\alpha {\mathcal P}_j(d)), \quad \mathcal{N}(d, \alpha) = \prod_{j=1}^{R} {\mathcal N}_j(d,\alpha).
\end{align}

Note that it follows from \cite[Lemma 1]{Radziwill&Sound} that the quantities defined in \eqref{defN} are all positive. Moreover, we deduce from \cite[Lemma 4.1]{Gao2021-2} that for any $\alpha \in \mr$,
\begin{align*}
%%\label{Nprodbound}
 \mathcal{N}(d, \alpha)\mathcal{N}(d, -\alpha)  \geq 1.
\end{align*}

  We apply the above to see that for two real numbers $c>0, k<0$,
\begin{align}
\label{Nbound}
\begin{split}
 & \sumstar_{(d,2)=1}{\mathcal N}(d, 2k) \Phi(\frac dX) \leq  \sumstar_{(d,2)=1}{\mathcal N}(d, 2k)
 {\mathcal N}(d, 2(k-1))^{c}{\mathcal N}(d, 2(1-k))^{c} \Phi(\frac dX) \\
\leq & \sumstar_{(d,2)=1}{\mathcal N}(d, 2k){\mathcal N}(d, 2(1-k))^{c/2}
 |L(\half, \chi_{8d})|^{c}{\mathcal N}(d, 2(k-1))^{c/2}  |L(\half, \chi_{8d})|^{-c} \Phi(\frac dX).
\end{split}
\end{align}

  For any fixed $k<0$, we choose $c$ to be any positive real number satisfying
\begin{align}
\label{cbound}
\begin{split}
 0<1-c(1-1/k)/2<1.
\end{split}
\end{align}

  We further note that (see \cite{R&Sound})
\begin{align*}
\begin{split}
   L(\half, \chi_{8d})=\sum_{n \geq 1}\frac {\chi_{8d}(n)}{\sqrt{n}} \in \mr.
\end{split}
\end{align*}
  It follows that we have $L(\half, \chi_{8d})^2=|L(\half, \chi_{8d})|^2$. Moreover, notice that $L(s, \chi_{8d})$ has an Euler product when $\Re s>1$ which implies that $L(s, \chi_{8d}) >0$ for any real $s>1$.  Thus, if one assumes RH, then one deduces by continuity that $L(s, \chi_{8d})\geq 0$ for any real $s \geq 1/2$.

 We now apply H\"older's inequality with exponents being $(1-c(1-1/k)/2)^{-1}, 2/c, -2k/c$ to the last expression in \eqref{Nbound} to see that
\begin{align}
\label{basiclowerbound}
\begin{split}
& \sumstar_{(d,2)=1}{\mathcal N}(d, 2k) \Phi(\frac dX) \\
 \leq & \Big ( \sumstar_{(d,2)=1} ({\mathcal N}(d, 2k){\mathcal N}(d, 2(1-k))^{c/2})^{(1-c(1-1/k)/2)^{-1}} \Phi(\frac dX)  \Big)^{1-c(1-1/k)/2} \\
& \times \Big ( \sumstar_{(d,2)=1} L(\half, \chi_{8d})^{2} {\mathcal N}(d, 2(k-1)) \Phi(\frac dX)  \Big)^{c/2}\Big ( \sumstar_{(d,2)=1} |L(\half, \chi_{8d})|^{2k}\Phi(\frac dX) \Big)^{-c/(2k)}.
\end{split}
\end{align}

  We point out here that for any $c$ satisfying \eqref{cbound}, we have
\begin{align*}
%%\label{basiclowerbound}
\begin{split}
& (2k+2(1-k)(c/2))(1-c(1-1/k)/2)^{-1}=2k.
\end{split}
\end{align*}

 Now, we deduce from \eqref{basiclowerbound} that in order to establish \eqref{orderofmagdag}, it suffices to prove the following two propositions.
\begin{proposition}
\label{Prop6}  Assume the truth of Chowla's Conjecture. Let $k$ be any negative real number and $X$ be a large real number. We have
\begin{align*}
%%\label{L2estmation}
\sumstar_{(d,2)=1} L(\frac 12,\chi_{8d})^2{\mathcal N}(d, 2k-2)\Phi\Big(\frac{d}{X}\Big)  \ll X ( \log X  )^{\frac {(2k)^2+2}{2}}.
\end{align*}
\end{proposition}

\begin{proposition}
\label{Prop5}  Assume the truth of Chowla's Conjecture. Let $k$ be any negative real number and $X$ be a large real number. We have for any positive $c$ satisfying \eqref{cbound},
\begin{align*}
%%\label{Nestmation}
\sumstar_{(d,2)=1} ({\mathcal N}(d, 2k){\mathcal N}(d, 2(1-k))^{c/2})^{(1-c(1-1/k)/2)^{-1}} \Phi\Big(\frac{d}{X}\Big) \ll & X ( \log X  )^{\frac {(2k)^2}{2}},  \\
\sumstar_{(d,2)=1}{\mathcal N}(d, 2k) \Phi(\frac dX)  \gg & X ( \log X  )^{\frac {(2k)^2}{2}}.
\end{align*}
\end{proposition}

  The above propositions in fact can be established following the proofs of Proposition 2.2-2.3 in \cite{Gao2021-3}. As the arguments are very similar, we omit them here.

  We note that the proof of \eqref{orderofmagdag} makes use of the observation that $L(\half, \chi_{8d})^2 \geq 0$. 
If one assumes RH, then our discussions above imply that $L(\half, \chi_{8d}) \geq 0$. We can then modify our arguments above to 
estimate quantities involving with $L(\half, \chi_{8d})$ only.  This can be further seen in our next proof of \eqref{orderofmagp}.

 For this, we fix a $k<0$ and choose $c'$ to be any positive real number satisfying
\begin{align}
\label{cbound1}
\begin{split}
 0<1-c'(1-1/(2k))<1.
\end{split}
\end{align}

  Then we have
\begin{align}
\label{Nbound1}
\begin{split}
 & \sum_{(p,2)=1}(\log p){\mathcal N}(p, 2k) \Phi(\frac pX) \leq  \sum_{(p,2)=1}(\log p){\mathcal N}(p, 2k)
 {\mathcal N}(p, 2k-1)^{c'}{\mathcal N}(p, 1-2k)^{c'} \Phi(\frac pX) \\
\leq & \sum_{(p,2)=1}{\mathcal N}(p, 2k){\mathcal N}(p, 1-2k)^{c'}
 |L(\half, \chi_{8p})|^{c'}{\mathcal N}(p, 2k-1)^{c'}  |L(\half, \chi_{8p})|^{-c'} \Phi(\frac pX).
\end{split}
\end{align}

  We apply H\"older's inequality with exponents being $(1-c'(1-1/(2k)))^{-1}, 1/c', -2k/c'$ to the last expression in \eqref{Nbound1} to see that
\begin{align}
\label{basiclowerbound1}
\begin{split}
& \sum_{(p,2)=1}(\log p){\mathcal N}(p, 2k) \Phi(\frac pX) \\
 \leq & \Big ( \sum_{(p,2)=1} (\log p)({\mathcal N}(p, 2k){\mathcal N}(p, 1-2k)^{c'})^{(1-c'(1-1/(2k)))^{-1}} \Phi(\frac pX)  \Big)^{1-c'(1-1/(2k))} \\
& \times \Big ( \sum_{(p,2)=1} (\log p) L(\half, \chi_{8p}) {\mathcal N}(p, 2k-1) \Phi(\frac pX)  \Big)^{c'}\Big ( \sum_{(p,2)=1} (\log p)|L(\half, \chi_{8p})|^{2k}\Phi(\frac pX) \Big)^{-c'/(2k)}.
\end{split}
\end{align}
  Here we note that from our discussions above, we have $L(\half, \chi_{8p}) \geq 0$ under RH.

  We also notice that for any $c'$ satisfying \eqref{cbound1}, we have
\begin{align*}
%%\label{basiclowerbound}
\begin{split}
& (2k+c'(1-2k))(1-c'(1-1/(2k)))^{-1}=2k.
\end{split}
\end{align*}

   Now, we deduce from \eqref{basiclowerbound1} that in order to establish \eqref{orderofmagp}, it suffices to prove the following two propositions.
\begin{proposition}
\label{Prop6p} Assume the truth of Chowla's Conjecture and assume the truth of RH for $\zeta(s)$ and GRH for $L(s, \chi_{8p})$ for all odd primes $p$,. Let $k$ be any negative real number and $X$ be a large real number. We have
\begin{align*}
%%\label{L2estmation}
\sum_{(p,2)=1} (\log p) L(\half, \chi_{8p}) {\mathcal N}(p, 2k-1) \Phi(\frac pX) \ll X ( \log X  )^{\frac {(2k)^2+1}{2}}.
\end{align*}
\end{proposition}

\begin{proposition}
\label{Prop5p} Assume the truth of Chowla's Conjecture and assume the truth of RH for $\zeta(s)$ and GRH for $L(s, \chi_{8p})$ for all odd primes $p$. Let $k$ be any negative real number and $X$ be a large real number. We have for any positive $c'$ satisfying \eqref{cbound1},
\begin{align*}
%%\label{Nestmation}
\sum_{(p,2)=1} (\log p)({\mathcal N}(p, 2k){\mathcal N}(p, 1-2k)^{c'})^{(1-c'(1-1/(2k)))^{-1}} \Phi(\frac pX)  \ll & X ( \log X  )^{\frac {(2k)^2}{2}},  \\
\sum_{(p,2)=1}(\log p){\mathcal N}(p, 2k) \Phi(\frac pX)  \gg & X ( \log X  )^{\frac {(2k)^2}{2}}.
\end{align*}
\end{proposition}

   The above propositions can be similarly established following the proofs of Proposition 3.4 and Proposition 3.6 in \cite{G&Zhao11} 
so that we again omit the arguments here.  This completes the proof of Theorem \ref{thmlowerbound}.

\vspace*{.5cm}

\noindent{\bf Acknowledgments.}  P. G. is supported in part by NSFC grant 11871082.

\bibliography{biblio}
\bibliographystyle{amsxport}

\vspace*{.5cm}

\end{document}